\documentclass{proc-l}

\usepackage{amsmath, amsthm, amssymb}

\usepackage{graphicx}

\usepackage[cmtip,all]{xy}

\usepackage{comment}

\newtheorem{theorem}{Theorem}[section]
\newtheorem{lemma}[theorem]{Lemma}
\newtheorem{prop}[theorem]{Proposition}
\newtheorem{Coro}[theorem]{Corollary}
\newtheorem{Conj}[theorem]{Conjecture}
\newtheorem{Prob}[theorem]{Problem}

\theoremstyle{definition}
\newtheorem{definition}[theorem]{Definition}
\newtheorem{example}[theorem]{Example}

\theoremstyle{remark}

\numberwithin{equation}{section}

\begin{document}

\title{Optimal presentations of Mathematical Objects %over Symbol Libraries
}


%    author one information
\author{Akshunna Shaurya Dogra}
\address{26 Vassar Street, Cambridge, MA-02142}
\curraddr{}
\email{adogra@mit.edu}
\thanks{}

%    \subjclass is required.
%\subjclass[2010]{Number Theory}

%\date{November 5, 2018}

%\dedicatory{}

%    "Communicated by" -- provide editor's name; required.
%\commby{}

%    Abstract is required.
\begin{abstract}
We discuss the optimal presentations of mathematical objects under well defined symbol libraries. We shall examine what light our chosen symbol libraries and syntax shed upon the objects they represent. A major part of this work will focus on discrete sets, particularly the natural numbers, with results that describe the presentation of the natural numbers under specific symbol libraries and what those presentations may reveal about the properties of the natural numbers themselves. We shall present bounds and constraints on the length and shape of presentations, connect already existing problems in other fields of mathematics to questions relevant to these presentations and otherwise illuminate why such a study can produce exciting results. 
\end{abstract}

\maketitle

\section{Introduction}
Let us describe a game. A middle schooler is pondering the principles of arithmetic. He has finally begun to understand what symbols like $1, 2, 3, +, *, \wedge$ represent. Suddenly, a flash of realization - the many symbols of the decimal system are superfluous, when it comes to writing the natural numbers. $1,+$ alone could do the job. However, it takes too long to write even the moderately large numbers with just $+$, so he decides to throw $*$ in the mix. Suddenly, even the largest numbers he can think of are easily writable, while still requiring very few symbols. Perhaps adding $\wedge$ would make them even more easier to write. How about that weird symbol $!$ that his teacher introduced in class recently? Soon, he discovers there are many ways of writing most numbers. In time, he asks himself a grand question: If the only symbols allowed to him were $1$ and some assortment of functions, operators and other tools from the mathematical world, what would be the most optimal way of writing a number $n \in \mathbb{N}$ for them?

A major inspiration for this work was a similar game played and later somewhat formalized by us at New York University \cite{A}. In section F26 of his famous book collecting standing problems of Number Theory, R. K. Guy gives a collection of problems emerging from another very similar game \cite{B}. Harry Altman made a couple of such games the central point of his studies in graduate school and wrote his thesis on them \cite{C}. Excepting different choices of counting and aesthetics, all were born from the same motivations expressed by the child above.

Like many other areas of mathematics, these simple games are actually ripe ground for seeding questions of paramount importance. In \cite{D}, J. Iraids et al. used computers to play the game at a level beyond humans and made observations that suggest a strong connection to the problem of the infinity of Sophie Germain primes. In \cite{E}, Blum Shub Smale connected a problem of this type to the P-NP problem. In their conclusions to \cite{F}, Zeilberger and Gnang express the opinion that the methods and results born from/used to tackle such problems could be powerful tools in Combinatorial, Computational and Experimental Mathematics. Guy and Conway in particular, seemed to have had strong interests in the questions that emerged from their explorations of this niche field, carrying on a torch one presumes was left to them by Selfridge, who was reportedly another enthusiast.

Many giants of Mathematics in the last century or so have found themselves in the vicinity of the question above. Yet, despite the intermittent interest expressed by several people, including some recent progress by Melo and Svaiter in 1996 \cite{G}, we have found scant evidence of any major successful effort to solve the standing problems in this field or even recognition of their importance. Indeed, our education in the matter has had to be piecemeal and patchy due to the absolute lack of any standard reference or direction.

We hope for the reader to get three things out of this work: \\
1. an understanding of the progress we have made in the field and make clear the connections this work has to the fragmented literature already existing, \\
2. forward a sensible, rigorous, yet diverse framework, from which future progress might be carried out on the many unsolved questions of the author and others,\\
3. entice the reader into doing research in an area that is young, naive and virtually unknown, while possibly offering keys to some of the most significant standing problems our times.

\section{Definitions:}
Mathematics requires a shared language to report the progress of ideas. Lack of clarity in that shared language is a hindrance to our progress. We focus on optimal presentations of mathematical objects under some defined symbol libraries and syntax. In an ideal world, we would define everything from the ground up, with no assumptions made as to the knowledge of the reader. However, even the greatest mathematicians have always worked first and made it rigorous later, guided by their intuition of the subject they were studying. Fortunately, it is certain that the targeted audience of this journal is more than capable of grasping the beginning steps, even when it requires a certain level of intuition and implicit knowledge and promises to build everything thereafter solely in terms of definitions already introduced. With that promise, we are ready to define the following:
\begin{definition}\label{Symlib}
A \textit{Symbol Library} $O$ is the collection of some symbols usable in a mathematical sense. $O$ is static, if its elements and composition are independent of its use. Else, $O$ is dynamic. $O$ is finite if it has only a finite number of elements.
\end{definition}
In a naive sense, a symbol library is a collection of symbols we can use to perform a mathematical operation or represent some mathematical object or a combination thereof. One should realize of course that the words $symbol$, $mathematical$ and $operation$ are ambiguous to a certain extent. Fortunately, context will clarify the coming discussions% and hence, for now, we shall skip the days we could devote to maybe coming up with a better definition%
. Symbol libraries in this work are aptly described in set notation. 
\begin{example}\label{ExSymlib}
Say the symbols allowed to us are $1$, $S$ $(successor)$, $+$, $*$. We can then write $O = \{1,S,+,*\}$. 
\end{example}
The most apparent use of a symbol library is creating objects by using the symbols in some arrangement or the other. The library of Example 2.2 can represent any $n \in \mathbb{N}$ in at least one way, by simply writing $\underbrace{SSS...SSS}_\text{$n-1$ copies}1$. Notice, this makes sense only under agreed upon use of syntax between reader and writer. For example $\underbrace{SSS...SSS}_\text{n copies}$ can be agreed upon to represent $n$ just as easily. However, we exclusively choose the former over latter because of convention and the fact that under prefix notation $*S1SS1$ is easily understood to be 6, whereas $*SSSSS$ could also be 4, without further clarification. Not all ways of writing out things are equally clear and useful. We shall adhere to mathematical convention and notation wherever existing and try to create the most intuitive and close conventions wherever we need to make our own. Let us define a term.
\begin{definition}\label{term}
A \textit{term} under $O$ is a string made using symbols exclusively in $O$.
\end{definition}
Individual terms made from a symbol library may or may not make sense to a reader. Symbols could be used in such an arbitrary manner, so as to throw this whole enterprise off course. Therefore, there is a distinction between some string that is simply a term and a string that represents something useful (and/or usefully). They shall be the guiding force of this work.
\begin{definition}\label{rep}
A \textit{presentation} under $O$, is a term evaluating to some mathematical object (using some agreed upon syntax).
\end{definition}
Under $O$ and syntax of Example 2.2, $*S1SS1$ is a presentation of 6, but $*SSSSS$ is a term. Notice $SSSSS1$, $+S1SSS1$, $+SS1SS1$ and $*1SSSSS1$ represent 6 too. Indeed, Gnang et al. in \cite{H} describe how the number of permissible presentations of $n \in \mathbb{N}$ grow as $n$ grows. This is not limited to just the natural numbers. In combinatorics, superpermutations of $n$ objects are strings that have all permutations of those $n$ objects as a substring.
\begin{example}\label{Exrep}
Let $O = \{1,2,3\}$ and consider all permutations of $1,2,3$. Both $123121321$ and $123231312121213132321$ are presentations of a superpermutation of $1,2,3$.
\end{example}
Even within presentations, some are more worthy of interest than others. In particular, we want to know something about the \textit{optimal} presentations of objects. However, before we go ahead and define them, we need to focus on the $optimal$ part of this discussion. The reasons are two-fold, one aesthetic, one functional:

1. there are presentation properties other than string length. In \cite{C} and \cite{D}, the presentations of interest are those using the least number of $1$ for $O = \{1, +, *\}$. In \cite{E}, the presentations of interest for some $n\in \mathbb{N}$ are those that minimize computation length of $n$ under $O = \{1, +, -, *\}$.

2. the symbols making up $O$ are usually different kinds of mathematical objects. In Example 2.2, $1$ is a natural number, while the other symbols are operations. Just like \cite{A} and \cite{C} decided upon different definitions of what an $optimal$ presentation of a number was, based on how they looked at the symbols available to them, so could others differ when working in a different setting.

In this work, we shall assume that the optimal presentation(s) for some object are the presentation(s) having the shortest string length (i.e. least number of symbols used), except when explicitly noted otherwise. So let us define optimal presentations in a general way, while keeping that convention at the back of our mind.
\begin{definition}\label{oprep}
An \textit{optimal presentation} of $n$ under $O$, is a presentation outputting the optimal numerical value for some pre-defined property of the syntax, amongst all permissible presentations of some mathematical object $n$. We call that optimized numerical value $k$, the $complexity$ of the element and write $c_O(n) = k$.
\end{definition}
The syntax property to be optimized can be quite arbitrary in general. However, in this work, we shall look at syntax properties whose values always $\in \mathbb{N}$. It is important to keep in mind that the syntax properties are a feature of $O$ and its usage, not of the object under study.

Let us say the mathematical object under study belongs to a set with some defined partial order or more. We might be interested in asking if the order of complexity is in any way related to the order of the elements. Notice that, \textit{a priori}, there is no reason to hope for such a relationship, as the order of the set under study might be defined in ways unrelated to whatever syntax property is asked to be optimized. However, a question that can be asked from a %discrete%
poset is: when does a particular complexity first appear? (that is, what is the minimal element(s) with a particular complexity \textit{k}). One may also be able to ask for the highest order element with complexity $k$. These elements can turn out to be extremely interesting. 
\begin{definition}\label{minugly}
\textit{Minimal element of k} is the minimal element(s) with complexity $k$ under some $O$. We denote it as $u_k$
\end{definition}

\begin{definition}\label{Max}
\textit{Maximal element of k} is the maximal element(s) which has a presentation with syntax property output $k$. We denote it as $M_k$
\end{definition}

Let us summarize the conventions we shall use hereon, before we start proving results. Unless otherwise specified, we assume that every term under $O$ can represent at most one unique element, that $S$ means a/the successor function(s) and the default syntax property output under consideration for a presentation is string length (number of symbols used). Finally, we shall assume that when we reference an object $a$ in a term, without writing out which presentation of $a$ we are using, we are writing $a$ in one of its optimal presentations by simply writing $a$. (For example, under $O$ of \ref{ExSymlib}, *(2)(3) means *S1SS1).
Let us prove some basic results while adhering to the conventions established.

\begin{prop}\label{Maxcomp_k}
Consider some set $A$ such that $\forall a\in A$, $Sa$ exists. If $S \in O$ and $M_k$ exists, then $c_O(M_k) = k$ and $u_k \leq M_k$, if $u_k$ exists in the same (sub)chain.
\begin{proof}
$c_O(M_k) \leq k$. If $c_O(M_k) = k-b$, $\underbrace{SS...SS}_\text{b copies}M_k > M_k$ is a presentation of length $k$, which implies $M_k$ is not a Maximal element. Hence, $c_O(M_k) = k$. (Note that $S$ is a stand in for all unique successor function(s) permissible in $A$. A unique successor function is taken to mean a function that produces at least one element that no other successor function can).
\end{proof}
\begin{Coro}\label{existu&M}
Let an infinite set $A$ have a finite number of unique successor functions and minimal elements. Let $O$ be finite, static and contain all minimal element(s) and successor function(s) of $A$. Then, $M_k$ exists $\forall k \in \mathbb{N}$.
%\begin{proof}
% M_k has to exist and as O is both finite and static and therefore so does u_k
%\end{proof}
\end{Coro}
\end{prop}
Consider the following examples detailing certain facts about $u_k,M_k$ under specially chosen $O$. Here, $P$ is the \textit{predecessor} function.
\begin{example}\label{specialcases}
Under $O = \{S, P, 1\}, M_2 = 2, u_2 = 0$, if we are looking at $\mathbb{Z}$.\\
Under $O = 2\mathbb{Z_-}\cup \{S\}, M_2 = -1$ and $u_2$ doesn't exist in $\mathbb{Z}$.\\
Under $O =\mathbb{Z_-}\cup\{+,*\}, M_3$ does not exist and $u_3=0$ in $\mathbb{Z}$.\\
Under any finite, static $O$ satisfying $\{1, S\}\subset O$, when making $\mathbb{N}, M_k$ exists $\forall k \in \mathbb{N}$. As shall be seen in \ref{noukskip}, so does $u_k$.
\end{example}

\begin{prop}\label{loglowerboun}
Under all static $O$ satisfying $|O| = c \in \mathbb{N}$, only $\frac{c(c^{k+1} - 1)}{c-1}$ elements can have complexity $\leq k$.
\begin{proof}
Simply sum up the number of possible terms of length $1,2,3,...k-1,k$. Notice, this implies $\geq\mathcal{O}(log)$ complexity upper bound on many kinds of sets under any kind of finite, static $O$. \end{proof}
\end{prop}
\begin{Coro}\label{statfinOukexist}
Let $O$ be finite and static. $\exists a \in A$, $c_O(a)=k'\implies \exists u_{k'}\in A$.
% only a finite number of k' length objects can exist due to O's finiteness. At least some of them have to be the minimal one with that length.
\end{Coro}

\ref{loglowerboun} stands irrespective of what operations are thrown into $O$, which is a little unexpected before we derive it, since one could potentially try to use exceedingly powerful operations to reduce complexity. But, while it does tell us that at least one element has higher than some log order complexity, it does not tell us much about how many. We will now prove a slightly less general, but more powerful result, by generalizing the arguments relayed to Gnang et al. in \cite{H}, by Noga Alon.
\begin{theorem}\label{loglowercount}
Let $A$ be a poset satisfying the following: $\forall n\in \mathbb{N}, \exists a \in A$ with order $n$.% and $\forall b\in A, \exists m\in \mathbb{N}$ such that $b$ has order $m$.%
(that is to say, $\forall n \in \mathbb{N}, \exists a \in A$ such that $a$ is reached by succeeding $n-1$ times from some minimal element $a_1 \in A$). Let $O$ be finite and static.\\% and all $b \in A$ are reachable by succeeding from some minimal element $b_0 \in A$, $m-1$ times for some $m\in \mathbb{N}$). \\
Then, for $\epsilon > 0$, number of elements $x\in A$, with order $\leq n$, such that $c_{O}(x) \leq (1-\epsilon)\log(n)/\log(|O|)$ is at most $\mathcal{O}(n^{1-\epsilon})$.  
\begin{proof}
%By standard assumptions, $O$ is finite and static $\implies O$ has finite number of minimal elements of $A$. 
Fix some $\epsilon>0$. Let $|O| = c, (1 - \epsilon)/\log(c) = z$. Let $o(x)$ be order of $x$ and\\
$E(n) = \{x |o(x) \leq n, c\textsubscript{O}(x) \leq z\log(n)\}$,\\
$k = max\{c\textsubscript{O}(x)| x \in E(n)\} \implies k \leq z\log(n)$.\\
$|E(n)| < \sum_{1 \leq i \leq k+1} c^{i} = \frac{c}{c-1}(c^{k+2}-1) = \mathcal{O}(n^{1-\epsilon})$ as $c^k \leq c^{zlog(n)} = n^{1-\epsilon}$\\
If even a single term is not a presentation, $k$ can replace $k+1$ in the sum.
\end{proof}
\end{theorem}
\begin{Coro}
Let $O$ be finite and static. Under any such $O$, the upper bound on the complexity of $\mathbb{N}$ is at least $\mathcal{O}(logn)$ .
\end{Coro}

\begin{prop}\label{arbitrary}
Let $A$ have only a finite number of unique chains $C_1, C_2,...., C_j$ and satisfy $A=C_1\cup C_2\cup ...C_j$ ($A$ is made of chains that are proper subsets of $A$, but of no other chain). Let $O$ contain all unique successor (predecessor) functions and minimal(maximal) elements of A.
\begin{enumerate}
\item If O contains only the successor(predecessor) functions and minimal(maximal) elements of $A,O$ is finite and static.
\item If complexity of $A$ is unbounded over $O$, there exists at least one unique chain in $A$, such that there exist arbitrarily large number of elements with order between successive elements of any complexity $k$ in that chain.
\item If in every possible unique chain of $A$, $\exists k\in \mathbb{N}$, such that the number of possible elements, with order between two successive elements of complexity $k$ is bounded, complexity of $A$ is bounded over $O$.
\end{enumerate}
\begin{proof}
Notice that number of total possible minimal(maximal) elements $\leq j$ and $A$ has to have a finite number of minimal elements. Similarly, $A$ can't have an infinite number of $unique$ successor functions (here a unique successor function is taken to mean a function that produces at least one element that no other successor function can). $(1)$ follows directly. Further, it implies each $a \in A$ has finitely many unique successor elements. \\
If all possible chains have some finite upper bounds on the length of gaps between consecutive elements of some complexity $k$, all chains can be totally covered with some finite number of successors applied to elements of this complexity $k$. Hence, all of $A$ is covered within some finite complexity. (2) and (3) follow.\\
\end{proof}
\end{prop}

\begin{theorem}\label{endsinS}
%Let $A$ satisfy: %$\forall a \in A, Sa$ exists and 
%all $a\in A$ have at least one immediately preceding element in $A$ (excepting minimal elements of $A$).
Let $O$ contain all unique successor function(s) and minimal element(s) of $A$. Then, for $k\geq 2$, %for all $k\in \mathbb{N}$, 
all $u_k\in A$ have at least one optimal presentation ending in a(the) successor function, if:
\begin{enumerate}
\item $\exists u_k\in A\implies$ there exists at least one element $u_k ^P\in A$ that is an immediate predecessor of $u_k$, with at least one presentation under $O$.
\item $\exists a \in A$, $c_O(a)=k'\implies \exists u_{k'} \in A$.
%\item at least one immediately preceding element of every $u_k$ has a presentation under $O$.
%\item $O$ is finite and static
%\item If %$O$ is static and infinite and 
%$\exists a \in A$, such that $c_O(a)=k$, then $\exists u_k \in A$.
\end{enumerate}
\end{theorem}
The existence and properties of $u_k$ in $A$ are of extreme importance in this work. Observe (2) need be assumed only if $O$ is infinite and/or not static, due to \ref{statfinOukexist}.

\ref{endsinS} shall require an important intermediate statement regarding the order of complexity of elements, relative to the $u_k$ elements.

\begin{prop}\label{MinUglythm}
%Let $u_k, u_{k'}$ exist. If $k<k', u_k < u_{k'}$
Let $u_{k'}$ exist for some $k' \in \mathbb{N}$. If $u<u_{k'}$ has at least one presentation under $O$, $c_O(u)<k'$. 
\begin{proof}
%Let us assume the theorem is not true and pick the smallest $u_{k'}$, such that $\exists u < u_{k'}$ and $c_O(u) = k > k'$. We have that $u_k<u$ exists. The complexity of any predecessor of $u_k$ has to be 
Let us assume there exist $k',k \in \mathbb{N}$ and $u\in A$, satisfying $k'<k$, $u <u_{k'}$ and $c_O(u)=k$. Therefore, for all such $k, \exists u_k\in A$ and so does at least one preceding element $u_k ^P$ for it. Pick the smallest such $u_k$. $c_O(u_k ^P) = c \geq k-1$ or we could use the appropriate successor function in $O$ to make $u_k$ with a smaller complexity. On the other hand, $c> k$ would imply $u_c\leq u_k ^P<u_k<u_{k'}$, which is a contradiction, as $k'<k<c$ and thus the $(k',c)$ pair also provides a counter-example. $c=k$ implies $u_k$ is not the minimal element of $k$. Therefore, $k' \leq c=k-1 < k$. $k'=c$ implies $u_{k'}$ is not the minimal element of $k'$. Therefore, $k'<c<k$ which again implies an element smaller than $u_k$ that provides a valid counter-example pair ($c,k'$). Hence, we have a contradiction and no such $k',k$ can exist. 
\end{proof}
\end{prop}
\begin{Coro}\label{beforeuk}
If $u_k$ exists, at least one preceding element of $u_k$ in $A$ under $O$ has complexity $k-1$.
% Su_k ^P = u_k which means the preceding element has exactly one less complexity
\end{Coro}
\ref{endsinS} follows from \ref{beforeuk}.
\begin{Coro}
There exist $u_1, u_2, ...u_{k-1},u_k\in A$, if $c_O(a) = k$, for some $a\in A$.
%\begin{proof}
% Pick u_k. u_k_P has complexity k-1 as O has all successors of $A$. Therefore u_{k-1} exists. Now pick u_{k-1_P}. Therefore u_k-2 exists and so on.
%\end{proof}
\end{Coro}

\begin{Coro}\label{noukskip}
$\forall k \in \mathbb{N}$, $\exists u_k\in A$, if the complexity of $A$ over $O$ is unbounded.
%\begin{proof}
%There are no holes in $u_k$ because $u_k -1$ has to have complexity $k-1$ for every $k$. Let us first observe that $\forall k\in \mathbb{N}, \exists M_k \in A \implies u_k$ exists by \ref{Maxcomp_k}
%\end{proof}
\end{Coro}
\section{Optimal presentations and $\mathbb{N}$}
The results presented until now are quite general in nature, assuming little about the cardinality, topology or other internal structure of the sets or symbol libraries in question. However, $\mathbb{N}$ was the original inspiration of this work and a lot of those results are generalizations of important facts about $\mathbb{N}$ under certain kinds of $O$. We will now explore the setting in which people before us have asked interesting and important questions. Hereon, we will almost always consider $\mathbb{N}$ with the caveat that $\{1,S\} \subseteq O$, a choice due to the following reasons:

1. Peano axioms have their roots in viewing $S$ as the generating function of $\mathbb{N}$, along with the principle of induction. Indeed, Hermann Grassmann first showed in \cite{I} that many facts of arithmetic could be revealed simply from a study of the successor function and its many implications on $\mathbb{N}$.

2. All of Hyperoperation theory and its related subjects find their origin in $S$, whether as a generator of higher order operations usable on some subset of $\mathbb{R}$ or otherwise.

3. The choice of $1$ is to a certain degree a matter of aesthetics. It could be replaced with any other number without major qualitative effects. However, $0$ or $1$ are the natural starting points, if we are studying \textit{complexity} and \textit{optimal} presentations of $\mathbb{N}$. We pick $1$ to ensure that every number is at most as complex as itself, no matter the $O$. The presence of $S$ places an important constraint locally: complexity can rise at most one at a time.

$\mathbb{N}$ under $\{1,S\}$ or $\{1,S,+\}$ is not a gainful setting. Interesting patterns first emerge when $*$ gets thrown into the mix. Indeed, a majority of the field as it exists now, has been a study of the following $O$: $\{1,S,*\}, \{1,+,*\}, \{1,S,+,*\}$ and $\{1,+,*,\wedge\}$. Hereon, the two major syntax property outputs will be:

1. number of $1$s and $S$s used in the string, used in \cite{B}, \cite{C} and \cite{D} (technically, results in \cite{B}, \cite{C} and \cite{D} count the number of $1$s used under $O=\{1,+,*\}$, but if we are seeking that output, adding an $S$ to $O$ and modifying our approach as above changes nothing).

2. length of a string (number of symbols used in total), used in \cite{A} and \cite{H}.\\
By default, we shall use output 2 and provide the results from naturally extending the arguments to 1. If we are giving a result using 1, we shall denote the symbol library as $O_|$, a nod to the beginnings of arithmetic, the use of $|$ marks to count. 

\subsection{$\mathbb{N}$ under arithmetical $O$ (and/or $O_|)$}
\begin{prop}\label{upper}
Let $\{1,S,+,*\} \subseteq O$. For all $n\geq 3, c_O(n) \leq 6log_3(n) - 3$ and $c_{O_|}(n) \leq 3log_2(n).$
\end{prop}
\begin{prop}\label{lower}
Let $O = \{1,S,+,*\}$. For all $n, c_O(n) \geq 5log_4(n) - 1$ and $c_{O_|}(n) \geq 3log_3(n).$
\end{prop}
The proofs of the two statements above will require us to establish some intermediate results interesting in their own right.
\begin{lemma}\label{u2u}
$u_{k+6} \geq 3u_k + 2$ and $u_{k+4} \geq 2u_k + 1$.  
\begin{proof}
Observe \ref{MinUglythm} ensures that $u_k \in \mathbb{N}$ are ordered by $k$. \\
This implies $c_O(n\leq u_k) \leq k$. Therefore, for all $n \leq 3u_k + 2$, one may write a presentation of the form $3a + r$ as $\underbrace{S...}_\text{\textit{r} copies}*SS1(a)$, where $1\leq a \leq u_k, 0\leq r \leq 2$. This tells us no number requiring complexity more than $k+6$ can appear until $3u_k+2$. Similar arguments give us the other inequality.
\end{proof}
\end{lemma}

The proof of \ref{upper} follows from \ref{u2u} applied recursively to the fact that $u_3\geq 3$. The proof above also tells us the following two facts: 
\begin{Coro}\label{Sophieprime}
$c_O(\left \lfloor{n/2}\right \rfloor) \geq c_O(n)-4$ (Alternatively, $c_{O_|}(\left \lfloor{n/2}\right \rfloor) \geq c_{O_|}(n)-3$).
\end{Coro}
\begin{Coro}
At most four consecutive natural numbers can be minimal elements for different values of $k$, for $k>3$. Further, if $u_k$ is odd, $u_{k+3}>u_k+3$, for $k>3$.
\end{Coro}
In \cite{A} and \cite{D}, the respective authors made the observation that $u_k$ under $O_|$ expectedly have a very strong tendency to be prime. More interesting was the observation in \cite{D} that the tendency of $\left \lfloor{u_k/2}\right \rfloor$ to be prime was about just as strong. \ref{Sophieprime} is interesting because it tells us that $\left \lfloor{u_k/2}\right \rfloor$ should be expected to be quite close to the minimal number of its complexity (the lower bound gains about complexity $2$ from $n$ to $2n$, so the upper bound on average must gain around the same). Hence, it has a complexity close to the maximum it could possibly have. There is no definitive work yet on the probabilities of $n$ being prime, based on how close to its respective $u_k$ it is, but it is an interesting direction to consider for anyone interested in the conjectured infinity of Sophie Germain primes.

\begin{prop}
If $O = \{1,S,+,*\}$, $M_k$ exists for all $k\in \mathbb{N}$, $c_O(M_k) = k$ and the following hold:\\
\begin{enumerate}
\item no subterms can have a S*()() form, i.e., no addition outside *, if * is present.
\item no subterms of form SSSSSS1 $(7)$ or numbers higher than $7$.
\item At most four subterms of form SS1 and at most one of S1, SSSS1 and SSSSS1.
\item S1 and SSSS1 can't be present in the same $M_k$.
\item S1 and SSSSS1 can't be present in the same $M_k$.
\item SSSS1 and SSSSS1 can't be present in the same $M_k$.
\end{enumerate}
\begin{proof}
$M_k$ has to exist by \ref{existu&M} and \ref{Maxcomp_k} gives it complexity $k$.\\
Property $1$ follows from $a(b+c),(a+c)b>ab+c$ if $a,b>1$ and $c\in\mathbb{N}$. \\
Property $1$ also implies that any $M_k$ has the form $*()*()*()....*()()$ where $()$ are stand ins for elements having the form $SS...SS1$. As $SSSSSS1$ is always replaceable by $*SS1SS1$ without change in term length, it can't be part of any $M_k$ or $M_k$ will not be the maximal element with a term length $k$. Similar arguments generate (4), (5), (6) and second half of (2). Lastly, $*SS1*SS1*SS1*SS1SS1$ and $*SSS1*SSS1*SSS1SSS1$ have the same length and $3^5 < 4^4$.
\end{proof}
\end{prop}

\begin{Coro}\label{formlower}
Let $k \geq 11$ and $m$ be the smallest integer such that $k\leq5m-1$. If $r = 5m - 1 - k$, $M_k = 3^r4^{m-r}$
\begin{proof}
To begin, pick $k\geq 50$. Notice only SSS1 can have an arbitrary amount of repetitions at this point and therefore at least 4 copies of it exist. Further, if any $S1$ subterms existed in $M_k$, we could replace $*S1SSS1$ by $*SS1SS1$ ($2*4$ by $3*3$) without change in term length. If any $SSSS1$ subterms existed, we could replace $*SSSS1*SSS1SSS1$ by $*SS1*SS1*SS1SS1$ ($5*4*4$ by $3*3*3*3$). Finally $*SSSSS1SSS1$ can be replaced by $*SS1*SS1SS1$ ($6*4$ by $3*3*3$). Therefore, $M_k$ has to have a form $3^a4^b$ for some $a,b \in \mathbb{N}$.\\
As noted before, five copies of $SS1$ are replaceable with four copies of $SSS1$. The only possible form left available to $M_k$, given that term length $\in \mathbb{N}$, is the one mentioned. $M_k$ exists, therefore it must have that form. A computer can verify the result for the remaining finite cases. 
\end{proof}
\end{Coro}

The proof of \ref{lower} for $k\geq 11$ follows from observing $4^c \geq 3^r4^{m-r}$, where $c=(k+1)/5$. $4^c \geq 3^r4^{m-r}$ is obtainable from introductory level applications of calculus. (Indeed, the same methods can also lead us to \ref{lower}, as was done in \cite{A}). The remaining finite cases can then be checked by hand or machine. The corresponding results for $O_| = \{1,S,+,*\}$ can be similarly obtained.\\

In \cite{C}, the \textit{defect} $\delta(n)$ of a number $n$, is defined as the difference between the actual and lowest possible complexity of a number $n$ under $\{1,+,*\}$. We borrow that notion and use \ref{formlower} and \ref{lower} together to give:
\begin{Coro}
The numbers with $0$ defect under $O = \{1,S,+,*\}$ are of the form $4^k$ (Alternatively, numbers with $0$ defect under $O_|$ are of the form $3^k$).
\end{Coro}
\begin{Coro}
Let $O=\{1,S,+,*\}$. For all primes $p>16, \delta (p)\geq 0.5$. Further, if $p>16$ and not of the form $2^k + 1$, $\delta (p)>1$. If $O=\{1,S,*\}, c_O(p) = c_O(p-1)$.
%basically to make a prime, even if you have +, you have to jump from some other number. Even if you're jumping from a defectless 4^k kinda number, the lower bund barely rises while complexity is 1 more. therefore, 0.5 if 4^k form and 1 if otherwise
\end{Coro}

A long standing conjecture in the field pertains to the long term behavior of \textit{complexity}, asking whether it tends to the lower bounds presented in \ref{lower}, as $n \to \infty$. 
\begin{Conj}\label{asymp}
$\frac{c_{O}(n)}{5log_4(n)} \to 1$, as $n \to \infty$ (Alternatively, $\frac{c_{O_|}(n)}{3log_3(n)} \to 1$, $n \to \infty$).
\end{Conj}

Results presented in \cite{D} suggest that the best upper bound for complexity might have a co-efficient around $10\%$ larger than the lower bound, based on asymptotic behavior predicted by computational and experimental examination. On the other hand, in \cite{A}, $\left \lceil{(5 + \frac{log4}{loga})log_4(n) + 1}\right \rceil$ was found to be a very good upper bound for the numbers studied, where $a\in \mathbb{N}$ was the smallest number such that $a^a \geq n$. $loga$ diverges, but extremely slowly, so \cite{A} and \cite{D} are at odds with each other regarding the prediction they make for the eventual fate of complexity, but don't present enough data to categorically rule out the other one's estimate.

Another long standing problem in the field asks if the complexity of $2^a$ is always $2a$, if we are counting the number of $1$s in the string to ascertain complexity. 
\begin{Conj}\label{3powk}
$\forall k \in \mathbb{N}, c_{O}(3^k) = 4k - 1$ (Alternatively, $\forall k \in \mathbb{N}, c_{O_|}(2^k) = 2k$).
\end{Conj}

\ref{lower} directly implies:
\begin{Coro}
Either \ref{asymp} or \ref{3powk} is false.
\end{Coro}
The problem of finding the best, generally applicable, upper bounds has stood the test of time and the results in \ref{upper} are the best known (Reyna and Lune in \cite{J} provided some much tighter upper bounds on complexity, but only for some subsets of density 1 in $\mathbb{N}$). We would now like to present a possible pathway to finding better upper bounds and link this subject to some questions in number theory. Consider the following conjecture. 
\begin{Conj}\label{Akshconjugly}
Let $O$ (or $O_|$) $= \{1,S,+,*\}$. For all $k \in \mathbb{N}, (u_k)^2 < u_{2k+\delta}$, for some fixed $\delta \in \mathbb{Z}$. In particular, $\forall k > K, (u_k)^2 < u_{2k+1}$, for some large enough $K$.
\end{Conj}
In our efforts to prove \ref{Akshconjugly} so that we could improve the upper bounds, we realized that $\delta$ did not necessarily have to be fixed for those purposes, merely \textit{small} compared to $k$ (say order $(logk)^a$ perhaps, for some fixed $a$). From those motivations rose the following question, which will lead us to make an interesting conjecture for the realms of number theory.

\begin{Prob}\label{produgly}
For $2 \leq k << N/2$, what is the length of the longest sequence of consecutive natural numbers $a_i\in \{1, 2, 3, ...N^2 -1, N^2\}$, such that for all $i, a_i$ does not satisfy: $a_i = xy$, ($x,y \leq kN$)
\end{Prob}
Such questions were made popular by Erdos and have been pursued by mathematicians for at least the last half century. Ford's Theorem 1 in \cite{P} gives us that the longest such sequence must have a length of order at least $(log N/log k)^\delta$ where $\delta = 1-\frac{(1+loglog 2)}{log 2}$.

However, the best upper bound we could find for the same problem was an order $\sqrt[]{N}$ bound we deduced by simply putting numbers of the form $M^2 - a^2$ between perfect squares from $1$ to $N^2$. Similar techniques might furnish a $N^{1/n}$ bound for $n\in \mathbb{N}$, but it is our contention that the best upper bound should be expected to be much smaller in order. Let us see how the problem above could connect with our work and how its resolution would be very beneficial for us.

The prime number theorem and its generalization by Landau \cite{L} allows us to count the number of unique numbers we may make with primes $<kN$. Count possible numbers of the form $p_1p_2$, where $p_1$ is a prime in $[(k-c-1)N,(k-c)N]$, where $0\leq c\leq k-1$ and $p_2$ is a prime $<N/(k-c)$. Then, start counting possible unique numbers of form $p_{21}p_2$ where $p_2$ is the same as before but $p_{21}$ is a semiprime in $[(k-c-1)N,(k-c)N]$. If we count out all the possible combinations of different multiplicities permissible, we get at least $(N/logN)^2$ unique numbers $<N^2$, that can satisfy the condition in \ref{produgly}. 

A problem similar in flavor is figuring out the length of the longest possible sequence of consecutive natural numbers $<N$, such that they are all composite. It is conjectured that the order is $(logN)^2$ \cite{K}. 
We conjecture the following:
\begin{Conj}\label{compconjcopy}
Let $a_i$ satisfy the following: it is a sequence of consecutive numbers $\in \{1, 2, 3, ...N^2 -1, N^2\}$, such that, $\forall i$, $a_i\neq xy$, for any $x,y<kN$.\\
Then, the length of $a_i$ is at most of order $(logN)^\delta$, for some fixed $\delta$.
\end{Conj}

\ref{compconjcopy} comes roughly from the same motivations that expect length of longest sequence of consecutive composite numbers $<N$ to be $(logN)^2$ and the shapes and complexity of optimal presentations (numbers). \ref{Akshconjugly} is motivated from computational analysis and data in \cite{A} and \cite{D}, in addition to \ref{upper}. 

If either one of \ref{compconjcopy} or \ref{Akshconjugly} is true, we would have the capacity to provide much tighter, general upper bounds. One may also ask if the coefficient $2$ in subscript in \ref{Akshconjugly} could be replaced by a smaller number (say $1.99$) as $n^2$ grows. The answer to the latter question is no and is related to the \ref{lower}. Finally, one could ask if the $\delta$ could take a negative value and what that would mean for the long term behavior of complexity (would that imply upper bounds converging to the lower bound?). We have been able to answer none of these questions definitively but would like to present whatever progress we have been able to make on them.

To begin, observe that the coefficient $2$ in \ref{Akshconjugly} need not be any higher. This is apparent from the observation that even if the upper bound was initially equal to the lower bound, if complexity rose fast enough from $n$ to $n^2$ to need a coefficient larger than $2$ to satisfy the inequality, we would end up violating the upper bound established in \ref{upper} at some point (indeed, we would violate any upper bound with $\mathcal{O}(log)$). Similarly, any $\mathcal{O}(log)$ upper bound would eventually end up lower than the lower bound in \ref{lower} if the co-efficient in the inequality in $\ref{Akshconjugly}$ was lower than $2$.

Assume \ref{Akshconjugly} is true. Therefore, if $c_O(u_{k_1}\leq n\leq u_{k_2}) \leq alog_4(n) - \delta$ for some $a \in \mathbb{R}$, $c_O(u_{2k_1+\delta}\leq m\leq u_{2k_2+\delta}) \leq alog_4(m) - \delta$. By starting at a large enough $n$ and iteratively progressing with $k_2=2k_1+\delta$, we ensure that all $m>u_{k_1}$ obey the upper bound (it would be in our interest to pick the best possible $a$ in the beginning step). As an application of the conjecture, we observe that the co-efficients in \ref{upper} would decrease to about $4.75$ and $2.30$, using the $u_k$ reported in \cite{A} and \cite{D}. 

Lastly, observe that if $\delta$ could be arbitrarily negative, it would still need to decrease slowly in value. One only need go beyond $n$ large enough, such that you can pick a small enough $a$, so that $alog(n)+|\delta|$ fits as an upper bound. In his thesis \cite{C}, Altman showed that the defect of $n$ is unbounded (if it were otherwise, the long term behavior of complexity would be already known). This tells us that there are limits on how negative $\delta$, even if \ref{Akshconjugly} turns out to be a weaker statement than the actual truth, vis a vis $\delta$. Consider this for example: $5log_4(n) + 10$ works as an upper bound until $n=4.5\times 10^6$, as investigated in \cite{A}. If $\delta$ dropped to a value of $10$ or lower before that, the upper bound mentioned above would stand for all $n$, contradicting the unboundedness of the defect. Hence, even if \ref{Akshconjugly} turns out to be weaker than the truth, vis a vis $\delta$, we would expect whatever function governs the $\delta$ values to not decrease too fast.

Now, assume \ref{compconjcopy} is true, with $\delta \geq 0$. Let $u_k$ be the smallest minimal element $>N$. Therefore, by \ref{u2u} and the argument carried out while discussing \ref{produgly}, at least $(N/logN)^2$ unique numbers smaller than $N^2$ have complexity $\leq 2k+3$ and the largest possible distance between two such successive numbers is order $(logN)^\delta$. Pick the smallest number $a$ such that $alog(m)$ is a valid complexity upper bound for $\sqrt[]{N}\leq m \leq N$. Therefore, $$c_O(N\leq n \leq N^2) < a[log(n)+ log(log(n))] + alog(\delta) + 5 + c$$ where $c$ is the smallest number such that $alog(m) + c$ is a valid upper bound for the first $\sqrt[]{N}$ numbers. Let $C = 5+alog(\delta)+c$. Now, pick the smallest $\epsilon$, such that $\epsilon log(n)>alog(log(n))$ for all $n$. For the same reasons governing the result due to \ref{Akshconjugly}, $(a+\epsilon)log(z)+C$ stands as a valid complexity upper bound $\forall z > N^2$.

The motivations for \ref{compconjcopy} and \ref{Akshconjugly} and the results they imply are telling us an important fact about the complexity of natural numbers: their complexity is self optimizing as we go to larger and larger $n$. This is partly because the number of valid presentations (or alternatively, paths) to a particular number increases rapidly as $n$ grows, as was shown in \cite{H} and partly because each of those presentations (or paths) use already optimized presentations (paths). 

Computationally obtaining the optimal presentations and complexity of $n$ requires the same knowledge only of $n-1$ and $a$ where $a$ is stand in for all possible factors of $n$. Hence, algorithms requiring runtime of only $\mathcal{O}(nlogn)$ easily compute and generate both optimal presentations and complexities till $n$.

In this section, we have focused on finite symbol libraries. A question that naturally arises is regarding the effects of infinite symbol libraries on the complexity of $\mathbb{N}$. For example, if $O = \{1,S,+\}\cup \mathbb{P}$, the complexity of $\mathbb{N}$ is bounded, by Helfgott's result on the ternary Goldbach conjecture \cite{O}(there are older results that would establish said fact too). However, if $O = \{1,S,*\}\cup \mathbb{P}$, the situation changes drastically.

We have seen examples of connections between problems in number theory and the optimal presentation of $\mathbb{N}$. We would like to present something in the reverse direction, a generalization of the arbitrariness of the prime gap using the methods shown in this work. We shall first need the following theorem:

\begin{theorem}\label{unboundcomP}
Let $O = \{1,S,*\}\cup \mathbb{P}$. Then, there exists $u_k \in \mathbb{N}$, for every $k\in \mathbb{N}$.
\end{theorem}
This is an interesting result, pointing out the fact that replacing $+$ by $*$ actually makes the complexity of $\mathbb{N}$ unbounded, even though $*$ is the more powerful operation by most standards. It also allows us a fascinating fact regarding the natural numbers:
\begin{Coro}\label{Multigapthm}
Let $k\in \mathbb{N}$ and $g_{m_k}$ be the difference between $m^{th}$ and $(m+1)^{th}$ number with multiplicity $k$ or lower. $g_{m_k}$ is unbounded for all $k \in \mathbb{N}$.

\end{Coro}
The corollary above is a generalization of the famous result that the prime gap can be arbitrarily large. It is obtainable from Landau's generalization of the Prime Number Theorem \cite{L}, but our route to it is shorter and far more elementary. To achieve that, we first provide the following lemmas using elementary combinatorial arguments, emulating the path laid out by Idris Mercer in \cite{M} and then applying on it a famous result first obtained by Euler:
\begin{lemma}\label{pnt}
Let $\pi (x)$ represent the prime counting function.\\
For all $x\geq 2, \pi (x) < 20x/log(x)$.
\begin{proof}
For all $n \in \mathbb{N}$, $\binom{2n}{n}$ $<4^n$ and $\binom{2n}{n}$ $\in \mathbb{N}$. \\
All prime numbers in $[n+1,2n]$ are factors of $\binom{2n}{n}$, as it is a natural number and the denominator does not have even a single prime from $[n+1,2n]$. Therefore,\\
$$n^{\pi (2n) - \pi (n)}< \binom{2n}{n} <4^n \implies \pi (2^{k}) - \pi (2^{k-1}) < \frac{2^k}{k-1}$$
Let us say we wanted to put a bound on the number of primes until some $k = 2m$, starting from $2$. The left part of the final inequality above gives us a telescopic sum that can be used to evaluate just what we want for our choice of $k$, giving: $$\pi (2^{2m}) - \pi (2) < 2^{m+1} + \frac{2^{2m+1}}{m} \implies \pi (4^m)< 1+2^{m+1}+\frac{2.4^m}{m}$$
Lastly, we use the fact that $\forall x \in \mathbb{N}$, we have $m\in \mathbb{N}$, such that $$4^{m-1}<x\leq 4^m \implies m-1<log_4(x)\leq m$$
Picking such an $m$ gives us:
$$\pi (x) < 1 +\frac{8x}{log_4(x)} + 4x^{1/2} < \frac{20x}{logx}$$
\end{proof}
\end{lemma}
\begin{lemma}\label{mnt}
Let $\pi _k (x)$ represent the k-multiplicity counting function, i.e., a function that counts numbers smaller than $x$ with multiplicity $k$.\\
For all $x\geq 2, \pi _k (x) < 20kx\frac{(loglog(x))^{k-1}}{log(x)}$.
\begin{proof}
We shall provide the proof for $k=2,3$ and trust the readers can carry out the induction thereafter themselves.
Let $p,q$ be primes such that $pq\leq n, p\leq q$. Therefore, $p\leq \sqrt[]{n}$. For any such $p$, we have $q \leq n/p$. Finally, by \ref{pnt} and $n/p\geq \sqrt[]{n}$, number of such primes is $<\frac{20n}{plog(n/p)}<\frac{40n}{plog(n)}$. To get all possible numbers, we sum over all possible $p<\sqrt[]{n}$. Let $P_n$ be set of all primes $\leq n$. We finally have:
$$\pi _2(n)<\sum _{p\in P_{\sqrt[]{n}}}\frac{40n}{plog(n)} = \frac{40n}{log(n)}\sum _{p\in P_{\sqrt[]{n}}}\frac{1}{p}$$
But $\sum _{p\in P_{\sqrt[]{n}}}\frac{1}{p} \sim loglog(n)$, by Euler's result on sum of reciprocal of primes (later more rigorously obtained by Mertens in \cite{N}). This identity gives us our result for $k=2$.

Let us now consider numbers of the form $pqr\leq n$, with $p\leq q\leq r$ being primes. Therefore, $p\leq n^{1/3}$. Further, $qr\leq n/p$ is a semi-prime. From the $k=2$ case, number of such semi-primes is $$\pi _2(n/p)< 40(n/p)(loglog(n)/log(n/p))<40(n/p)(loglog(n)/log(n^{2/3}))$$
Therefore 
$$\pi _3(n/p)< 40(n)(3/2)(loglog(n)/log(n))\sum _{p\in P_{\sqrt[3]{n}}}\frac{1}{p}$$
Summing over all possible $p$ gives us the result for $k=3$ case. This process can then be repeated iteratively to get to any required $k\in \mathbb{N}$.
\end{proof}
\end{lemma}

We now have all the tools needed to prove \ref{unboundcomP}.
\begin{proof}
Let us assume \ref{unboundcomP} is false. It could be false in two ways: $u_k$ does not exist for some types of $k$ or the complexity of $\mathbb{N}$ over $O$ is bounded above by some $K\in \mathbb{N}$.
\begin{enumerate}
\item \ref{noukskip} ensures that %every existing $u_k$ has an optimal presentation ending in $S$. Therefore, $c_O(u_k - 1) = k-1$, if $u_k$ exists, which implies $u_{k-1}$ exists. Thus, 
if $c_O(\mathbb{N})$ is unbounded, then $u_k$ exists for all $k \in \mathbb{N}$.\\
\item Let $c_O(\mathbb{N}) < K$ for some smallest $K\in \mathbb{N}$. Every $n\in \mathbb{N}$ has at least one optimal presentation under $O$. Let us start counting the number of total possible optimal presentations of any length until some large $n \in \mathbb{N}$.
\begin{itemize}
\item Length 1: $\pi (n)+1$: the only numbers with optimal presentations of this length are the primes and $1$.
\item Length 2: $<\pi (n)$: the only numbers with possible optimal presentations of this length are numbers succeeding primes.
\item Length 3: $<\pi (n) + \pi _2(n)$: the only numbers with possible optimal presentations of this length are numbers two more than primes and semiprimes.
\item Length 4: $<\pi (n) + \pi _2(n) + \pi _2(n): $ the only numbers with possible optimal presentations of this length are numbers three more than a prime or numbers one more than a semiprime or numbers that are products of a prime and a number succeeding a prime.
\item Length 5: ............................\\
.............................................\\
.............................................\\
.............................................\\
.............................................\\
\item Length $K-1$: $<\pi (n)+\pi _2 (n)+...................+\pi _{\lfloor K/2 \rfloor}(n)+\pi _{\lceil K/2 \rceil}(n)$ : the only numbers with possible optimal presentations of this length are numbers $K-2$ more than a prime or numbers $K-4$ more than a semi prime or ................. or numbers with multiplicity $\lfloor K/2 \rfloor$ (and numbers with multiplicity $\lceil K/2 \rceil$, if $K$ is odd.)

\end{itemize}
\end{enumerate}
By \ref{mnt}, the density of all possibilities in $n$ is smaller than: $$\frac{1+ \pi (n) + \pi (n) + \pi _2 (n)........+ \pi _{\lceil K/2 \rceil}(n)}{n} \to 0, \textnormal{     as } n\to \infty$$ 
But there are no other possibilities for optimal presentation forms left and the number of possible optimal presentations has to be larger than $n$ (hence, have density $>1$). We are left with a contradiction and hence no such $K$ can exist. 
\end{proof}

\ref{Multigapthm} follows from the following observation: if $\exists k,l\in \mathbb{N}$, such that the gap between consecutive numbers $m,n\in \mathbb{N}$ with multiplicity $\leq k$, is always $<l$, we could use $S\in O$, along with $\mathbb{P}$ to cover all of $\mathbb{N}$, with complexity at most $2k-1+l$, thus contradicting \ref{unboundcomP}. \\

\ref{unboundcomP} naturally extends even further:
\begin{Coro}
Let $O = \{1,S,*\}\cup \mathbb{P}_m$, where $\mathbb{P}_m$ is the set of all numbers with multiplicity $m$ or less. Then, there exists $u_k \in \mathbb{N}$, for every $k\in \mathbb{N}$.
\end{Coro}
It might interest the reader to observe that statements \ref{endsinS} and beyond from section $2$ apply to the system we just studied.

\textbf{% THE EVERY U_K ENDS IN S THEOREM SHOULD HOLD TRUE FOR O = 1,S,* U |P TOO. PROVE IT.
}

\bibliographystyle{amsplain}

\end{document}